\documentclass[12pt]{article}

\usepackage{enumerate}
\usepackage{amsthm,amsmath,amssymb}
\usepackage{graphicx}
\usepackage[colorlinks=true,citecolor=black,linkcolor=black,urlcolor=blue]{hyperref}
\usepackage[utf8]{inputenc}

\theoremstyle{plain}

\theoremstyle{definition}

\theoremstyle{remark}

\begin{document}

\title{ A Variation on the Homological Nerve Theorem }
\author{Luis Montejano }

\maketitle 

\begin{abstract}  An equivalent but useful version of the Homological Nerve Theorem is proved.
\end{abstract}

\section{ Introduction}

Let $X$ be a polyhedron,  $F=\{A_1,\dots A_m\}$  a polyhedral cover of $X$ and let $N=N(F)$ be the nerve of the family $F$. Denote by $N^{(k)}$ the k-skeleton of the simplicial complex $N$.
During these paper we are going to use reduced homology with coefficients in a field. We say that $A\subset X$ is $\rho$-acyclic  if $\tilde H_*(A)=0$ for $*\leq \rho$. Furthermore, $\tilde H_{-1}(A)=0$ means that $A$ is not empty.

\smallskip

The Homological Nerve Theorem, as stated by Meshulam in [3],  claims the following:

\smallskip

 \noindent {\bf Homological Nerve Theorem}. Suppose that for every $\sigma \in N^{(k)}$,  

\begin {itemize}

\item $\bigcap _\sigma$ is $(k-\mid\sigma \mid +1)$-acyclic.  \quad Then,

\end {itemize}

\begin{center}

Rank $\tilde H_{k+1}(N) \leq$ Rank $\tilde H_{k+1}(X).$

\end {center}

\noindent and for every $0\le j\le k,$

\begin{center}

$\tilde H_j(N) =\tilde H_j(X).$

\end {center}

The purpose of this paper is to prove the following equivalent but useful variation of the Homological Nerve Theorem:

\bigskip

\noindent {\bf Theorem 1.}  Suppose that for every $\sigma \in N^{(k)}$,  

\begin {itemize}

\item $\tilde H_{k-\mid\sigma \mid +1}(\bigcap _\sigma)=0$. Then,

\end {itemize}

\begin{enumerate}

\item Rank $\tilde H_{k+1}(N) \leq$ Rank $\tilde H_{k+1}(X),$
\item Rank $\tilde H_{k}(X) \leq$ Rank $\tilde H_k{(N)}.$

\end {enumerate}

\medskip

In many applications of the Homological Nerve Theorem the useful conclusion is that $\tilde H_{k+1}(X)=0$   implies $\tilde H_{k+1}(N)=0 $. So Theorem 1 helps to improve these 
 results since the hypothesis needed to achieve this conclusion are much weaker.  In Section 3, we give a couple of examples of this fact. 
 
 \section{ Killing the homology groups }

Let $X$ be a polyhedron and let $f:S^d\to X$ be a PL-map. Denote by $X\cup_f B^{d+1}$ the polyhedron obtained be attaching the $(d+1)$-cell to $X$. Then,  if $u\in \tilde H_d(S^d)$ is the fundamental class, and $f_*(u)\in \tilde H_d(X)$ is not zero, we have that for every $*\not=d$,
\begin {center}
$\tilde H_*(X)=\tilde H_*(X\cup_f B^{d+1})$, and

$\tilde H_d(X\cup_f B^{d+1})=\frac{\tilde H_d(X)}{<f_*(u>}$.
\end {center}
\medskip
 
 The next lemma is essential in several proofs of this paper. Its proof is part of the "folklore" and we include it by completeness.

\medskip

\noindent {\bf Lemma 1.}.   Let $A$ be a polyhedron. Then by attaching $\lambda$-cells to $A$, $\lambda\le k$, we obtain a polyhedron $\tilde A$ containing $A$ such that 
\begin{itemize}
\item $\tilde A$ is $(k-1)$-acyclic, and 
\item $\tilde H_i(A)=\tilde H_i(\tilde A)$, for $i\geq k.$
\end {itemize}
\medskip 
\noindent {\bf Proof. } If $k=1$, the theorem is trivial because we just have to connect components with arcs. If $k=2$, Let $g_1, \dots ,g_t:S^1\to A$ continuous maps such that $(g_i)_*(u)\not= 0$ and $\{(g_i)_*(u)\}_1^t$ generates $\tilde H_1(A)$, where $u\in \tilde H_1(S^1)$ is the fundamental class.  Then, by attaching 2-cell to $A$, via $g_1, \dots ,g_t:S^1\to A$ we achieve our purpose. Suppose now $k=3$. We first kill the fundamental group. 
Let $\psi_1, \dots ,\psi_t:(S^1, 1)\to (A,*)$ be continuous maps such that  $\{(\psi_i)_*(u)\}_1^t$ generates $\pi_1(A,*)$, where now $u\in \pi_1(S^1,*)$ is the generator. Then, by attaching 2-cells to $A$, via $\psi_1, \dots ,\psi_t:S^1\to A$ we obtain $\tilde A$,  killing the fundamental group of $A$ but perhaps creating $2$-dimensional homology. In other words $\pi_1(\tilde A,*)=0=\tilde H_1(\tilde A)$,  and $\tilde H_*(A)=\tilde H_*(\tilde A)$, for $*\geq 3$.  Now, by the Hurewitz Theorem $\pi_2(\tilde A,*)=\tilde H_2(\tilde A; Z)$. So let $\psi_1, \dots ,\psi_t:(S^2, 1)\to (\tilde A,*)$ be continuous maps such that   $\{(\psi_i)_*(u)\}_1^t$ generates $\pi_2(\tilde A,*)$ and $(\psi_i)_*(u)\not= 0$, where this time $u\in \pi_2(S^.*)$ is the generator. Then, by attaching 2-cells to $A$, via $\psi_1, \dots ,\psi_t:S^2\to A$ we obtain $\tilde {\tilde A}$.  First of all,  note that $\tilde H_2(\tilde {\tilde A}, Z)=0$ but also that, again by the Hurewicz Theorem, $(\psi_i)_*(u)\not= 0 \in \tilde H_2(\tilde {\tilde A};Z)$. Therefore, by the universal coefficient Theorem for homology, since $\tilde H_1(\tilde {\tilde A})=0$, we have that $\tilde H_2(\tilde {\tilde A})=0$. Furthermore, again by the universal coefficient Theorem,  $(\psi_i)_*(u)\not= 0 \in \tilde H_2(\tilde {\tilde A})$ and consequently $\tilde H_*( \tilde {\tilde A})=\tilde H_*(A)$, for $*\geq 3$.  The proof of the theorem for higher k's is completely analogous to the case $k=3$. \qed

\bigskip

\noindent {\bf Proof of Theorem 1.}  

The idea is to carefully kill the homology of the $A_i$'s  and its intersections by attaching cells in such a way that we can use the Homological Nerve Theorem.  We will do it in such a way that we do not modify the nerve of the family and in such a way we do not change the $(k+1)$-dimensional homology of the resulting new $X$.

\smallskip
\noindent {\it Inductive Claim $C_r$}.  $1\le r \le k+1$.

It is possible to construct polyhedra $\tilde A_i$,   and define $\tilde F=\{\tilde A_i\}_1^m$ and $\tilde X=\bigcup_1^m \tilde A_i $, in such a way that:
\begin{enumerate}

\item $ A_i\subset \tilde A_i$,
\item $N=N(\tilde F)$,
\item $\tilde H_{k+1}(\tilde X)=\tilde H_{k+1}(X)$,
\item Rank $\tilde H_{k}(X) \leq$ Rank $\tilde H_{k}(\tilde X),$
\item  for every $\sigma \in N^{(k)}$,  $\tilde H_{k-\mid\sigma \mid +1}(\bigcap_{i\in\sigma} \tilde A_i)=0,$
\item if $\sigma \in N^{(k)}$ and $\mid \sigma \mid \geq r$,  $\bigcap_{i\in\sigma} \tilde A_i$ is $(k-\mid\sigma \mid +1)$-acyclic.
\end {enumerate}

By defining $ A_i = \tilde A_i$, $C_{k+1}$ is clearly true.  Next we shall prove that if $C_{r+1}$ is true then $C_r$ is also true.  Note that if $C_1$ is true, then we obtain our result by applying the Homological Nerve Theorem, because on one side  
$\tilde H_{k+1}(\tilde X)=\tilde H_{k+1}(X)$, $\tilde H_{k+1}(N)=\tilde H_{k+1}(N(\tilde F))$ and Rank $\tilde H_{k+1}(N(\tilde F)) \leq$ Rank $\tilde H_{k+1}(\tilde X)$ and in the other side,  Rank $H_k (X) \leq$ Rank $\tilde H_{k}(\tilde X) =$ Rank $\tilde H_{k}(N),$

Suppose the inductive claim $C_{r+1}$ is true.  Let us first fix $\sigma \in N^{(k)}$ with $\mid \sigma \mid = r$. 
 Hence, by (4), $\tilde H_{k-r+1}(\bigcap_{i\in\sigma} \tilde A_i)=0$. By Lemma 1, we can kill of the $(k-r)$-homology of $\bigcap_{i\in\sigma} \tilde A_i$, by attaching $\lambda$-cells of dimension smaller or equal than $k-r+1$. 
So, for $i\in \sigma$, we obtain $\tilde {\tilde A}_{i}$ from $\tilde A_{i}$, by attaching the same $\lambda$-cells of dimension smaller or equal than $k-r+1$,  in such a way that  $\bigcap_{i\in \sigma} \tilde {\tilde A}_i$ is $(k-r+1)$-acyclic.
Finally,  for $\rho \in \{1, \dots , m\} -\sigma$, let $\tilde {\tilde A}_\rho=\tilde A_\rho$ and let $\tilde {\tilde F}=\{\tilde {\tilde A}_i\}_1^m$ and $\tilde {\tilde X}=\bigcup_1^m \tilde {\tilde A}_i$.

Our next purpose is to prove that for $S\subset \{1, \dots , m\}$,
$$\bigcap_{i\in S}\tilde A_i= \bigcap_{i\in S}\tilde {\tilde A}_i.$$
whenever $S$ is not contained in $\sigma$, and $ \bigcap_{i\in S}\tilde {\tilde A}_i$ is obtained from $ \bigcap_{i\in S}\tilde  A_i$ by attaching $\lambda$-cells of dimension smaller or equal than $k-r+1$, whenever $S\subset \sigma$.

Suppose first $S$ is not contained in $\sigma$, then 
$\bigcap_{i\in S}\tilde A_i= \bigcap_{i\in S}\tilde {\tilde A}_i,$
 because if $j\in S-\sigma$, then $\tilde {\tilde A}_j \cap \tilde X=\tilde {\tilde A}_j$.  Then  $\bigcap_{i\in S}\tilde {\tilde A}_i=  (\bigcap_{i\in S-\{j\}}\tilde {\tilde A}_i)\cap \tilde {\tilde A}_j$ $= (\bigcap_{i\in S-\{j\}}\tilde {\tilde A}_i)\cap (\tilde {\tilde A}_j\cap \tilde X)$=$(\bigcap_{i\in S}\tilde {\tilde A}_i)\cap \tilde X= \bigcap_{i\in S}(\tilde {\tilde A}_i\cap\tilde X)=\bigcap_{i\in S}\tilde A_i$.

On the other hand, if $S\subset \sigma$, then by definition of the $\tilde {\tilde A}_i$'s,  we have that $ \bigcap_{i\in S}\tilde {\tilde A}_i$ is obtained from $ \bigcap_{i\in S}\tilde  A_i$ by attaching $\lambda$-cells of dimension smaller or equal than $k-r+1$. 
\smallskip

So, here are some important  consequences of the above:
\begin {itemize}
\item $N(\tilde {\tilde F})=N$.

\item $\tilde H_{k+1}(\tilde {\tilde X})=\tilde H_{k+1}(X).$

\item Rank $\tilde H_{k}(X) \leq$ Rank $\tilde H_{k}(\tilde {\tilde X}).$

\item For every $\tau \in N^{(k)}$ with $\mid\tau\mid \geq r+1$, we have that 
$\bigcap_{i\in\tau}\tilde {\tilde A}_i=\bigcap_{i\in\tau}\tilde  A_i$.

\item For every $\tau \in N^{(k)}$ with $\mid\tau\mid = r$ and $\tau \not= \sigma$, we have that 
$\bigcap_{i\in\tau}\tilde {\tilde A}_i=\bigcap_{i\in\tau}\tilde  A_i$.

\item  For every $\tau \in N^{(k)}$ with $\mid\tau\mid < r$, we have that $\tilde H_{k-\mid\tau \mid +1}(\bigcap _{i\in\tau} \tilde {\tilde A}_i )=0$.  This is so, because either $\bigcap_{i\in\tau}\tilde {\tilde A}_i=\bigcap_{i\in\tau}\tilde  A_i$ or  $\bigcap _{i\in\tau} \tilde {\tilde A}_i $ is obtained from $\bigcap _{i\in\tau} \tilde A_i$ by attaching $\lambda$-cells of dimension $k-r+1 < k-\mid\tau \mid+1$.

\end{itemize}

By performing one by one this construction, to every $\sigma \in N^{(k)}$ with $\mid\sigma\mid=r$, we obtain that $C_r$ is true. This completes the proof of the theorem. \qed

\bigskip

Note that if for every $\sigma \in N^{(k)}$,  

\begin {itemize}

\item $\bigcap _\sigma$ is $(k-\mid\sigma \mid +1)$-acyclic.

\end {itemize}

Then, for every $0\le k^\prime \le k$, the following is true:

\noindent for every $\tau \in N^{(k^\prime)}$,  
$\tilde H_{k^\prime-\mid\tau \mid +1}(\bigcap _\sigma)=0$. This implies, by Theorem 1, that for every $0\le k^\prime \le k,$

\begin{enumerate}

\item Rank $\tilde H_{k^\prime+1}(N) \leq$ Rank $\tilde H_{k^\prime +1}(X),$
\item Rank $\tilde H_{k^\prime }(X) \leq$ Rank $\tilde H_{k^\prime}{(N)},$

\end {enumerate}

\noindent and therefore that  for every $0\le k^\prime \le k,$

\begin{center}

$\tilde H_{k^\prime}(N) =\tilde H_{k^\prime}(X)$  \quad and 

Rank $\tilde H_{k+1}(N) \leq$ Rank $\tilde H_{k+1}(X),$

\end {center}

\noindent thus proving that Theorem 1 is equivalent to the Homological Nerve Theorem.

\section{ Some consequences of Theorem 1}

In many applications of the Homological Nerve Theorem the useful conclusion is that $\tilde H_{k+1}(X)=0$   implies $\tilde H_{k+1}(N)=0 $. The purpose of this section is to show  a couple of examples in which 
the use of Theorem 1 instead of the Homological Nerve Theorem  helps to improve these results
 since the hypothesis needed to achieve the conclusions are much weaker. 

Let us start proving that the Topological Helly Theorem, obtained by Kallai and Meshulam in [2], can be derived from the Homological Nerve Theorem.
\bigskip

\noindent {\bf Topological Helly Theorem}  Let $F=\{A_1,\dots A_m\}$ be a collection of polyhedra in $ R^d$,  $m\geq d+2$. Suppose that for every subfamily $F^\prime$ of $F$ of size $n$, \quad $1\leq n\leq d+1$,

\begin {center}

 $\bigcap_{F^\prime}$ is $(d-n)$-acyclic,

\end {center}

\noindent then

$$ \bigcap_{F} \not= \emptyset.$$

\bigskip

\noindent {\bf Proof.}  As usual for a proof of a Helly type theorem, the proof follows by induction, and in this case also by a simple Mayer Vietories argument, from the case $m=d+2$.

Suppose $m=d+2$ and suppose also $ \bigcap_{F} = \emptyset.$  Then $N(F)$ is the boundary of a simplex with $d+2$ vertices and hence homeomorphic to the $d$-sphere.  On the other hand, if $k=d-1$, for every $\sigma \in N^{(d-1)}$,  
$\bigcap _\sigma$ is $(d-\mid\sigma \mid)$-acyclic.  Then, by the first conclusion of the Homological Nerve theorem,  $\tilde H_d(N(F))=0$ because $\tilde H_d(\bigcup _F)=0$, but this is a contradiction to the fact that $N(F)$ is a $d$-sphere. \qed

As we can see, we only used the first conclusion of the Homological Nerve Theorem, so exactly the same proof but now using Theorem 1 instead of the Homological Nerve Theorem yields the following topological Helly-type Theorem, first proved in [5].
\bigskip

\noindent {\bf Theorem 2}  Let $F=\{A_1,\dots A_m\}$ be a collection of polyhedra in $R^d$,  $m\geq d+2$. Suppose that for every subfamily $F^\prime$ of $F$ of size $n$, \quad $1\leq n\leq d+1$,

\begin {center}

 $\tilde H_{d-n}(\bigcap_{F^\prime})=0,$

\end {center}

\noindent then

$$ \bigcap_{F} \not= \emptyset.$$

Let $K$ be a simplicial complex. Suppose the vertices of $K$ are painted with $I=\{1, ... ,m\}$ colours, that is, there is a partition of the set of  vertices of $K$; $V(K)=V_1 \sqcup...\sqcup V_m$. A simplex $\sigma =\{v_1, ... , v_m\}\subset V(K)$ is rainbow if it contains exactly one vertex of every colour.  Finally, let $S\subset I$ be a subset of  colours.  Let $V_S\subset V(K)$ be the set of vertices of K painted with a colour in $S$ and let $K_S$ be the subcomplex of $K$ generated by vertices of $V_S$.

\medskip

The following theorem was proved by Meshulam [4] and also Aharoni-Berger [1]. 

\bigskip 

\noindent {\bf Theorem 3}. Let $K$ be a simplicial complex and suppose the vertices of $K$ are painted with $I=\{1, ... ,m\}$ colours. Then $K$ contains a rainbow simplex provided
\begin {center}
$K_S$ is $(s-2)$-acyclic,
\end{center}

\noindent for every subset $S\subset I$ of $s$ colours,  $1\leq s\leq m$.

\bigskip

In the proof of Theorem 3, given Meshulam [4], he used the first conclusion of the Homological Nerve Theorem to conclude the existence of a rainbow simplex, so 
exactly the same proof but now using Theorem 1 instead of the Homological Nerve Theorem yield the following improvement of Theorem 3.  For more about this kind of Sperner-type Theorems, see [6].
\medskip

\noindent {\bf Theorem 4}. Let $K$ be a simplicial complex and suppose the vertices of $K$ are painted with $I=\{1, ... ,m\}$ colours. Then $K$ contains a rainbow simplex provided
\begin {center}
$\tilde H_{s-2}(K_S)=0,$ 
\end{center}

\noindent for every subset $S\subset I$ of $s$ colours,  $1\leq s\leq m$.

\section{Acknowledgements}
The author wish to acknowledge  support  form CONACyT under 
project 166306 and  support from PAPIIT-UNAM under project IN112614.

\vspace{1cm}

\noindent \textsc{Luis Montejano }

\bigskip

\noindent \emph{Instituto de Matem\'{a}ticas, Unidad Juriquilla.\newline
National University of M\'{e}xico.}

\noindent E-mail address: \emph{luismontej@gmail.com}


\bigskip






\begin{thebibliography}{11}


\bibitem{AB} Aharoni R. and Berger E. The intersection of a matroid and a simplicial complex. Trans. Amer. Math. Soc., 358. No.~11, (2006), 4895-4917.

\bibitem{KM}  Kalai G. and Meshulam R. A topological colourfull Helly theorem. Adv. Math. 191. No.2. (2005), 305-311.

\bibitem{Mh1}  Meshulam R. The Clique Complex and Hypergraph Matching. Combinatorica. 21. No.1. (2001), 89-94.

\bibitem{Mh2}  Meshulam R.  Domination numbers and homology. Journal of Combinatorial Theory, Series A. 102. (2003), 321-330.

\bibitem{M} Montejano L. A new topological Helly theorem and some transversal results.  Discrete and Computational Geometry. 52. No.2. (2014), 390-398.

\bibitem{M} Montejano L. Homological Sperner-type Theorems.  Preprint.  2016

\end{thebibliography}
\end{document}